\documentclass[]{siamart220329}
\usepackage{amsmath}
\usepackage{amssymb}
\usepackage{algorithm}
\usepackage{algorithmic}
\usepackage{subfig}
\usepackage{graphicx}
\usepackage{mdframed}
\usepackage{mathtools}
\usepackage{float}
\usepackage{multirow}
\usepackage{matlab-prettifier}
\usepackage{algorithm}
\usepackage{algorithmic}

\newsiamremark{remark}{Remark}
\newsiamremark{hypothesis}{Hypothesis}
\crefname{hypothesis}{Hypothesis}{Hypotheses}
\newsiamthm{claim}{Claim}

\usepackage{graphicx} 
\usepackage{amsfonts,amssymb,amsmath}
\usepackage{cite} 
\usepackage{color}
\usepackage[font=small,labelfont=bf]{caption}

\usepackage{algorithm}
\usepackage{algorithmic}

\DeclareMathAlphabet{\mathbf}{OT1}{cmr}{bx}{it}

\newcommand{\vb}{\mathbf b}

\newcommand{\ve}{\mathbf e}
\newcommand{\vf}{\mathbf f}

\newcommand{\vr}{\mathbf r}
\newcommand{\vs}{\mathbf s}
\newcommand{\vt}{\mathbf t}
\newcommand{\vu}{\mathbf u}
\newcommand{\vv}{\mathbf v}
\newcommand{\vw}{\mathbf w}
\newcommand{\vx}{\mathbf x}
\newcommand{\vy}{\mathbf y}

\newcommand{\vnull}{\boldsymbol{0}}

\newcommand{\vV}{\mathbf V}
\newcommand{\vW}{\mathbf W}

\renewcommand{\d}{\,\mathrm d}

\DeclareMathOperator{\spn}{span}
\DeclareMathOperator{\sign}{sign}

\title{Krylov Subspace Recycling With Randomized Sketching For Matrix Functions
}
\author{Liam Burke\thanks{School Of Mathematics, Trinity College Dublin, College Green, Dublin 2, Ireland, \email{burkel8@tcd.ie}}
\and Stefan G\"{u}ttel\thanks{Department of Mathematics, The University of Manchester, M13 9PL Manchester, United Kingdom, \email{stefan.guettel@manchester.ac.uk}}.  S.\,G. acknowledges a fellowship from the Alan Turing Institute under the EPSRC grant EP/W001381/1 and a Royal Society Industry Fellowship IF/R1/231032.}

\usepackage{amsopn}

\ifpdf
\hypersetup{
  pdftitle={Krylov Subspace Recycling With Randomized Sketching For Matrix Functions},
  pdfauthor={Liam Burke and Stefan G\"{u}ttel}
}
\fi

\begin{document}

\maketitle

\begin{abstract}
A Krylov subspace recycling method for the efficient evaluation of a sequence of matrix functions acting  on a set of vectors is developed. The method improves over the recycling methods presented in [Burke et al., arXiv:2209.14163, 2022] in that it uses a closed-form expression for the augmented FOM approximants and hence circumvents the use of numerical quadrature. We further extend our method to use randomized sketching in order to avoid the arithmetic cost of orthogonalizing a full Krylov basis, offering an attractive solution to the fact that {recycling algorithms built from shifted augmented FOM} cannot easily be restarted. The efficacy of the proposed algorithms is demonstrated with numerical experiments.
\end{abstract}


\overfullrule=0pt 

\section{Introduction}
This paper is concerned with the development of Krylov subspace recycling algorithms for
the efficient computation of a sequence of matrix function applications of the form
\begin{equation}\label{eq:sequence_of_matrix_function_applications}
    f(A^{(i)}) \vb^{(i)}, \qquad i=1,2,\ldots
\end{equation}
with matrices $A^{(i)} \in \mathbb{C}^{N \times N}$ and vectors $\vb^{(i)} \in \mathbb{C}^{N}$.
This is a common problem arising in a variety of scientific computing applications. Common examples include simulations of quantum chromodynamics (QCD) which often require the evaluation of the matrix sign function $f(z) = \mbox{sign}(z)$ on a sequence of slowly changing Dirac matrices \cite{knechtli2017lattice, frommer2000numerical}, and the time integration of stiff systems of ordinary differential equations (ODEs) using exponential integrators, which  requires evaluation of the matrix exponential $f(z) = \mbox{exp}(z)$ on a sequence of matrices and vectors which may change only slightly between time steps~\cite{exp_integrators, Lu2003}.
Perhaps the most well-known example occurs when $f(z) = z^{-1}$, in which case (\ref{eq:sequence_of_matrix_function_applications}) is equivalent to solving a sequence of linear systems of equations; see e.g., \cite{parks2006recycling} and references therein.

In the special case where the matrices $A^{(i)}=A$ remain fixed, the number of different vectors $\vb^{(i)}$ is relatively small and they are all available simultaneously, {block Krylov subspace methods}~\cite{frommer2017block, frommer2020block} may be used for evaluating \eqref{eq:sequence_of_matrix_function_applications}. The recycling algorithms we develop in this work are of particular interest in the case when many vectors $\vb^{(i)}$ are available in sequence rather than simultaneously and/or the matrices~$A^{(i)}$ are slowly changing. The idea is to update an augmentation space from one problem to the next in order to improve convergence.  We stress that recycling is different from deflated restarting~\cite{Morgan2002}, the latter assuming that the matrix $A=A^{(i)}$ and vector $\vb=\vb^{(i)}$ are fixed while updating an augmentation space after each restart. 

The \emph{recycled FOM for functions of matrices} presented in \cite{burke2022krylov} (and therein abbreviated as r(FOM)$^2$) is the first Krylov subspace method to treat a sequence of matrix function applications of the form (\ref{eq:sequence_of_matrix_function_applications}) using subspace recycling. The method is based on augmented FOM approximants to the solution of the shifted linear systems appearing in an integral representation of~$f(A)\vb$, and it  uses numerical quadrature for the approximation of the integral. For notational convenience we refer to this algorithm as \emph{rFOM (quad)} throughout this paper. 
Unlike standard shifted FOM,  the residuals of augmented shifted FOM approximants are not necessarily collinear when the shift varies, making it difficult to restart rFOM (quad). The rFOM (quad) algorithm thus suffers from excessive storage and orthogonalization costs when the matrices $A^{(i)}$ are non-Hermitian. This is in addition to the need of deriving and evaluating an accurate quadrature formula for the integral representation of~$f(A)\vb$.

The aim of this paper is to overcome the limitations of quadrature-based recycling, developing a robust method that effectively makes use of the fact that Krylov information may be reused across the different problems \eqref{eq:sequence_of_matrix_function_applications} to speed up convergence.

\begin{figure}
    \centering
\includegraphics[width=.8\textwidth]{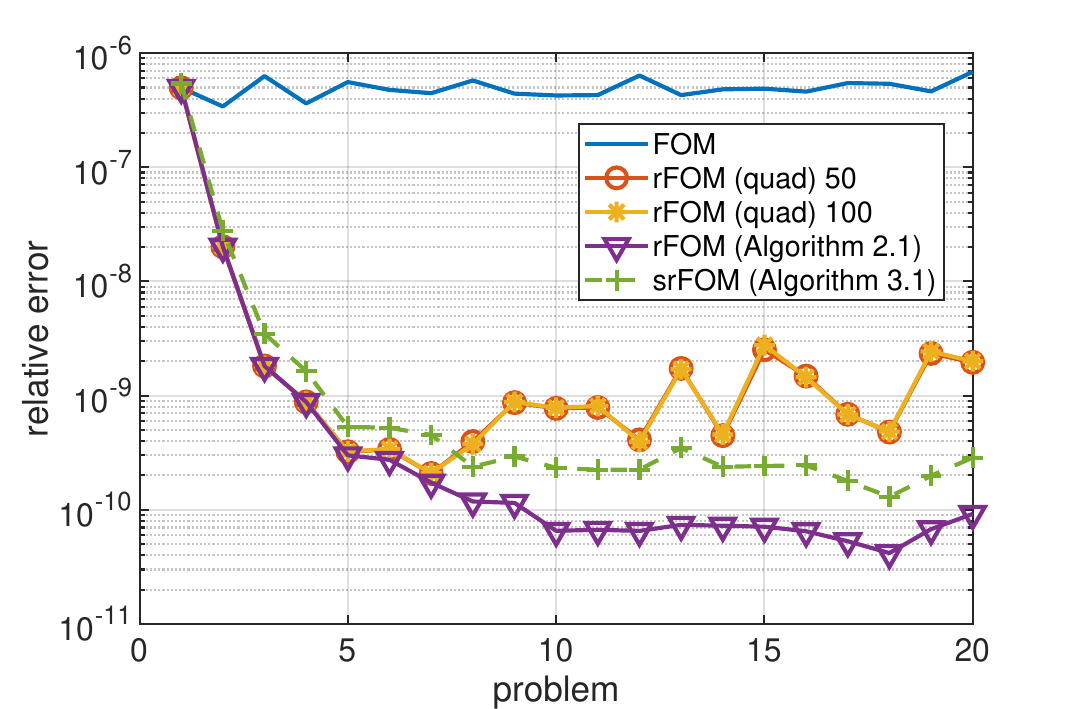} 
\caption{Solving $20$ QCD problems $(A^{(i)})^{-1/2}\vb^{(i)}$  with slowly changing matrices and varying vectors. We show the relative error  obtained for each problem with standard FOM (without recycling), the rFOM (quad) algorithm of \cite{burke2022krylov}, and the new algorithms rFOM and srFOM. The Arnoldi cycle length is $m = 100$ for all problems, and the recycling methods use a recycling subspace of  dimension of $k = 30$. The rFOM (quad) algorithm is run using $50$ and $100$ quadrature nodes, yielding the same convergence curves. This indicates that no  improvement can be achieved by further increasing the number of quadrature nodes. Algorithm~\ref{alg:rFOM} uses a closed form of the recycled FOM approximant and does not require quadrature. Algorithm~\ref{alg:srFOM} additionally utilizes randomized sketching to reduce orthogonalization cost and further speed up the computation at the expense of a slightly increased error.
} 
\label{fig:compare_quad_v_closed}%
\end{figure}

\smallskip
The overall structure and key contributions  are as follows. 
\smallskip
\begin{itemize}
    \item In section~\ref{sec:subs} we derive new closed-form expressions for augmented subspace approximations of matrix functions. The resulting method, referred to simply as rFOM,  is mathematically equivalent to the integral-based augmented FOM approximants in~\cite{burke2022krylov}, but it does not require  quadrature for its implementation and is observed to be numerically more robust for larger Arnoldi cycle lengths (see Figure~\ref{fig:compare_quad_v_closed}).  We also prove a convergence result for  rFOM  when applied for Stieltjes functions of Hermitian positive definite matrices. 
    \item In section~\ref{sec:sketch} we show how to combine rFOM with randomized sketching (resulting in srFOM), allowing for the use of nonorthogonal bases of the Krylov and augmentation spaces, thereby reducing computational complexity. We show how to utilize the sketched Rayleigh--Ritz approach from \cite{NakatsukasaTropp21} to update the {recycling subspace}.
    \item We also state a closed-form expression for an augmented GMRES-type approximant in section~\ref{sec:sketch} and its sketched counterpart. In comparison to the sGMRES approximant proposed in \cite{GS23}, our closed-form GMRES-type approximant does not require numerical quadrature for its evaluation. (However, both approximants are \emph{not} mathematically equivalent.)
    \item Section~\ref{sec:implem} discusses various implementation details of our methods, in particular, a SVD-based stabilization approach that demonstrably improves the numerical robustness of srFOM, and an error estimator that can be used to dynamically control the dimension of the Krylov space used for each problem.
    \item Finally, section~\ref{sec:numex} contains numerical experiments to demonstrate the efficacy of the new rFOM and srFOM  in comparison to the state-of-the-art. 
\end{itemize}

\smallskip

For a visual illustration of the algorithmic main contributions in this paper we refer to Figure~\ref{fig:compare_quad_v_closed} and its caption. The 20~QCD test problems solved in this plot are described in more detail in section~\ref{sec:qcd}.

\section{Generalized subspace extraction for matrix functions}\label{sec:subs}

The {rFOM (quad)} algorithm of \cite{burke2022krylov} is derived from  augmented  FOM approximants $\vx_m(\sigma)$ for the shifted linear systems 
\[
    (\sigma I - A) \vx(\sigma) = \vb
\]
appearing in the integral representation
\begin{equation}\label{eq:integral_representation}
    f(A) \vb = \int_{\Gamma}  (\sigma I - A)^{-1} \vb  \d \mu(\sigma) = \int_{\Gamma} \vx(\sigma)  \d \mu(\sigma)
\end{equation}
with some (possibly complex-valued) measure  $\d \mu$ defined on a contour $\Gamma\subset \mathbb{C}$. 
The augmented approximation~$\vf_{m}\approx f(A)\vb$ is obtained by replacing $\vx(\sigma)\to \vx_m(\sigma)$ in the integral representation and evaluating the integral using numerical quadrature. 
Here we show that it is not necessary to use quadrature, and that an augmented FOM algorithm and many other variants of subspace extraction can be implemented by using a closed-form expression involving the function~$f$.

Let the columns of two matrices $\textbf{V}_m$ and $\textbf{W}_m$ span a search and constraint space in $\mathbb{C}^N$, respectively. Then, following \cite[eq.~(5.5)]{burke2022krylov} with some notational simplifications, we define the  \emph{Bubnov--Galerkin approximant for the solution of $(\sigma I - A)\vx(\sigma) = \vb$} as 
\begin{equation}
\vx_m(\sigma) := \vV_m \left[  \vW_m^* (\sigma I - A) \vV_m \right]^{-1} \vW_m^* \vb.\label{eq:afom}
\end{equation}
Clearly, $\vx_m(\sigma) \in \spn(\vV_m)$ and  it is easy to verify that the residual 
\[
\vr_m(\sigma) := \vb - (\sigma I - A) \vx_m(\sigma)
\]
is orthogonal to $\spn(\vW_m)$:
\begin{eqnarray*}
\vW_m^* \vr_m(\sigma) &=& \vW_m^* [\vb - (\sigma I - A) \vx_m(\sigma) ]\\
 &=& \vW_m^* [\vb - (\sigma I - A)  \vV_m \left[  \vW_m^* (\sigma I - A) \vV_m \right]^{-1} \vW_m^* \vb ]\\
  &=& \vW_m^* \vb - \left[\vW_m^*  (\sigma I - A)  \vV_m\right] \left[  \vW_m^* (\sigma I - A) \vV_m \right]^{-1} \vW_m^* \vb \\
  &=& \vnull.
\end{eqnarray*}
We also refer to \cite{HH05} for an in-depth discussion of such (oblique) projection approaches for matrix functions.

We stress that the approximant \eqref{eq:afom} is very general. In particular, when $\vV_m$ spans the Krylov space $\mathcal{K}_m(A,\vb)=\spn \{\vb,A\vb,\ldots,A^{m-1}\vb\}$ and $\vV_m = \vW_m$, then   \eqref{eq:afom} is the approximant produced by the \emph{full orthorgonalization method} (FOM)~\cite{Saad2003}. If~$\vV_m$ spans the  Krylov space $\mathcal{K}_m(A,\vb)$ and $\vW_m = (\sigma I - A) \vV_m$, then \eqref{eq:afom} is the \emph{generalized minimal residual} (GMRES) approximant~\cite{SaadSchultz1986}. This is easy to verify by writing $\vx_m(\sigma)=\vV_m \vy_m(\sigma)$ and then checking that the condition $\vW_m^* \vr_m(\sigma)=\vnull$ is equivalent to the normal equations for $\|\vb - (\sigma I - A) \vV_m \vy_m(\sigma)\|_2\to\min_{\vy_m(\sigma)}$. 
However, it is usually impractical to have a constraint space $\vW_m$ that depends on the shift $\sigma$, so we also define a \emph{GMRES-type approximant} by setting $\vW_m = A \vV_m$ independent of~$\sigma$.  

\smallskip

Assuming that $\vW_m^*\vV_m$ is nonsingular, \eqref{eq:afom} can be rewritten as
\begin{eqnarray*}
\vx_m(\sigma) &=& \vV_m \left[  \sigma \vW_m^*\vV_{m} - \vW_m^* A \vV_m \right]^{-1} \vW_m^* \vb \\ 
&=& \vV_m \left[ \left( \sigma I - \vW_m^* A \vV_{m} (\vW_m^*\vV_m)^{-1} \right) \vW_m^*\vV_{m} \right]^{-1} \vW_m^* \vb \\ 
&=& \vV_m (\vW_m^*\vV_{m})^{-1} \left[ \sigma I - \vW_m^* A \vV_{m} (\vW_m^*\vV_m)^{-1} \right]^{-1} \vW_m^* \vb.
\end{eqnarray*}
Alternatively, and equivalently,
\begin{eqnarray}\label{eq:second_solution_approximation}
\vx_m(\sigma) &=& \vV_m  \left[ \sigma I - (\vW_m^*\vV_m)^{-1} \vW_m^* A \vV_m  \right]^{-1} (\vW_m^*\vV_m)^{-1} \vW_m^* \vb,
\end{eqnarray}
and we prefer to use \eqref{eq:second_solution_approximation} for now. 
Using the integral representation \eqref{eq:integral_representation}
and replacing $\vx(\sigma)$ by the approximant $\vx_m(\sigma)$  defined in \eqref{eq:afom}, 
we obtain 
\begin{equation}
\vf_m := \vV_m  f\left( (\vW_m^*\vV_m)^{-1}  \vW_m^* A \vV_m   \right)  (\vW_m^*\vV_m)^{-1} \vW_m^* \vb.
\label{eq:fm}
\end{equation}
We refer to $\vf_m$ as the \emph{Bubnov--Galerkin approximant of $f(A)\vb$}; see also~\cite{HH05}. When reading this formula, it might be helpful to keep in mind that $\vV_m (\vW_m^*\vV_m)^{-1}  \vW_m^*$ is the matrix representation of the oblique projector onto the span of $\vV_m$ along $\vW_m$.
It is remarkable that \eqref{eq:fm} is a closed formula and its evaluation does not require quadrature as opposed to the integral representation used in~\cite{burke2022krylov}.

Another attractive property of \eqref{eq:fm} is that it lends itself to randomized sketching, and we will exploit this in section~\ref{sec:sketch} below. Indeed, a similar closed-form representation has been derived for the sketched FOM approximant in \cite{GS23}, but \eqref{eq:fm} is more general as it also encompasses GMRES and GMRES-type approximants of $f(A)\vb$.

\subsection{Closed-form augmented FOM approximation}

We now assume that $\vV_m = \vW_m$ with linearly independent (but not necessarily orthonormal) columns. We then have
\[
(\vW_m^*\vV_m)^{-1} \vW_m^* = (\vV_m^*\vV_m)^{-1} \vV_m^* = \vV_m^\dagger,
\]
and therefore \eqref{eq:fm} simplifies to
\begin{equation}
\vf_m = \vV_m  f\left(\vV_m^\dagger A \vV_m   \right)  \vV_m^\dagger \vb.
\label{eq:afom2}
\end{equation}
If the columns of $\vV_m$ span the Krylov space $\mathcal{K}_m(A,\vb)$, this is the Rayleigh--Ritz representation of the FOM approximant as it has been studied extensively in the literature. In the context of subspace recycling, however,   $\vV_m$ may span an arbitrary augmentation space which is not necessarily a Krylov space. 

The practical evaluation of \eqref{eq:afom2} requires the computation of $\vV_m^\dagger A \vV_m$. This can either be done by using sketching (as we will discuss in section~\ref{sec:sketch}), or alternatively by orthogonalization of $\vV_m$.  Following the augmented FOM case studied in \cite{burke2022krylov}, $\vV_m = \vW_m = [ U, V_m]$ contains an augmentation basis $U\in\mathbb{C}^{N\times k}$ (such as $k$ approximants to a few relevant eigenvectors of $A$) and $V_m\in\mathbb{C}^{N\times m}$ is a Krylov basis. After orthogonalization $\vV_m \leftarrow \texttt{qr\_econ}(\vV_m)$ we can simply evaluate  \eqref{eq:afom2} as $\vf_m = \vV_m  f\left(\vV_m^* A \vV_m   \right)  \vV_m^* \vb$.

Our closed-form  recycled FOM is given in the Algorithm~\ref{alg:recyl} below, and Figure~\ref{fig:compare_quad_v_closed} demonstrates its numerical behavior on a sequence of QCD test problems. We will discuss this problem in greater detail in section~\ref{sec:numex}, but the key message is that closed-form FOM with recycling can yield significant convergence acceleration for sequences of matrix function computations, and it can also be more stable than the quadrature-based rFOM~(quad)~\cite{burke2022krylov}. We  believe that the improved stability of Algorithm~\ref{alg:recyl} results, at least partially, from the explicit orthogonalization of the augmented basis~$\vV_m$, which is in contrast to recycling variants that split $\vx_m(\sigma) = \vs_m(\sigma) + \vt_m(\sigma)$ and keep $U$ and $V_m$ separate, and also different from variants which work with $\vV_{m}$ but do not orthogonalize its columns.

\begin{remark}
In cases where the dimension $m$ is known a priori or kept fixed, it can be beneficial to  flip the order of $V_m$ and $U$, using $\vV_m = \vW_m = [ V_m, U]$. This is because the augmentation matrix $U$ typically has fewer columns than $V_m$ (i.e., $k \ll m$) and if the matrix $V_m$ has been computed by the Arnoldi process we already have an Arnoldi decomposition 
\begin{equation}\label{eq:arn}
A V_m = V_m H_m + h_{m+1,m} \vv_{m+1} \ve_m^T
\end{equation}
with $V_m^* V_m = I_m$, $V_m^* \vv_{m+1} = 0$, and an upper-Hessenberg matrix $H_m\in\mathbb{C}^{m\times m}$. Hence ensuring that $\vV_m = [V_m, U]$ has fully orthonormal columns requires the orthonormalization of only $k$ additional vectors. Further, for orthonormal $\vV_m$ we have
\begin{equation}
 \vV_m^* A \vV_m = 
\begin{bmatrix}
V_m^* A V_m & V_m^* A U \\
U^* A V_m & U^* A U 
\end{bmatrix}
=
\begin{bmatrix}
H_m & V_m^* A U \\
h_{m+1,m} (U^* \vv_{m+1}) \ve_m^T  & U^* A U 
\end{bmatrix}.
\label{eq:VAV}
\end{equation}
Hence, when care is taken to keep the Krylov basis numerically orthonormal (e.g., using the modified Gram--Schmidt process with reorthogonalization),  it is possible to reduce the number of matrix-vector products with $A$ in computing  $\vV_m^* A \vV_m$. 

In the special case when the matrices in the problem sequence \eqref{eq:sequence_of_matrix_function_applications} remain fixed, it is possible to update an Arnoldi-like decomposition for $A U$ and further reduce the number of matrix-vector products in forming~\eqref{eq:VAV}. This is the case, for example, when the new augmentation subspace $U^{(i+1)}$ is  obtained as
\[
U^{(i+1)} := \vV_m^{(i)} X = [ V_m^{(i)}, U^{(i)}] X
\]
from a previous cycle with some $X\in\mathbb{C}^{(m+k)\times k}$. If the product $Y^{(i)} := AU^{(i)}$ is available from the previous cycle in which have also computed an Arnoldi decomposition $A V_m^{(i)} = V_m^{(i)} H_m^{(i)} + h_{m+1,m}^{(i)} \vv_{m+1}^{(i)} \ve_m^T$, then
\[
A U^{(i+1)} = \left[ V_m^{(i)} (H_m^{(i)}X) + h_{m+1,m}^{(i)} \vv_{m+1}^{(i)} \ve_m^T X \, , \, Y^{(i)}X \right],
\]
which requires no additional matrix-vector products with~$A$. 
\end{remark}

\begin{algorithm}[t]
\caption{Closed-form recycled FOM for  $f({A}^{(i)})\vb^{(i)}$ \hfill rFOM \label{alg:rFOM}}
\begin{algorithmic}[1]
\STATE{\textbf{Input:} ${A}^{(i)}\in\mathbb{C}^{N\times N}$, $\vb^{(i)}\in\mathbb{C}^{N}$, integers $m$ and $k$}
\STATE{\textbf{Output:} Approximants $\widehat{ \vf}_{m}^{(i)} \approx f({A}^{(i)}) \vb^{(i)}$}\\[1mm]
\STATE{Initialize $U\leftarrow [\ ]$}
\FOR{$i = 1,2,\ldots$}
\STATE{Generate decomposition $A^{(i)} V_m = V_m H_m + h_{m+1,m} \vv_{m+1} \ve_m^T$ for $\mathcal{K}_{m}({A}^{(i)},\vb^{(i)})$}
\STATE{Set $\vV_m \leftarrow [ U, V_m]$}
\STATE{Orthonormalize $\vV_{m}$ and construct $\vV_m^* A^{(i)} \vV_m$}
\STATE{$\widehat{\vf}_{m}^{(i)} \leftarrow \vV_{m}  f(\vV_m^* A^{(i)} \vV_m)  \vV_m^* \vb^{(i)}$}
\smallskip
\STATE{Update $U$ by extracting $k$-dimensional subspace from $\spn(\vV_m)$}
\ENDFOR
\smallskip
\end{algorithmic}
\label{alg:recyl}
\end{algorithm}

\subsection{Convergence of augmented FOM} For the case that $A$ is a Hermitian positive definite matrix and $f$ obeys the Stieltjes integral representation
\[
    f(z) = \int_{0}^{+\infty} (z + t)^{-1} \d\mu(t)
\]
we can essentially follow the arguments in \cite{FrommerGuettelSchweitzer2014b} to prove convergence of the augmented FOM approximant $\vf_m$ defined in \eqref{eq:afom2} to $f(A)\vb$. We spell out the details below.

Let us assume that $\vV_m = [ V_m, U]$ where $V_m$ spans $\mathcal{K}_m(A,\vb)$ and $U$ spans the augmentation space. 
Let us introduce  $\vx(t) := (A + t I)^{-1}\vb$, $t\geq 0$, and the associated augmented FOM approximant $$\vx_m(t): =  \vV_m (\vV_m^\dagger A \vV_m + t I)^{-1}\vV_m^\dagger\vb$$ 
as well as the non-augmented (standard) FOM approximant
$$\underline\vx_m(t) :=  V_m (V_m^\dagger A V_m + t I)^{-1}V_m^\dagger\vb.\ \ $$ 
(Note the bold $\vV_m$ versus non-bold $V_m$ typesetting.) 
Associated with $\vx_m(t)$ are the residuals $\vr_m(t) := \vb - (A + tI)\vx_m(t)$ and errors $\ve_m(t) := \vx(t) - \vx_m(t)$. The vectors $\underline\vr_m(t)$ and $\underline\ve_m(t)$ are defined analogously.

Denoting by $\lambda_{\max}$ the largest eigenvalue of $A$ and defining the norm $\|\ve\|_{A} := \sqrt{\vv^* A \vv}$, it is easily verified that for $t\geq 0$ we have
\[
\|\vv\|_{A} \leq \sqrt{\lambda_{\max}/(\lambda_{\max} + t)}\cdot \|\vv\|_{A + tI}.
\]
Further, the orthogonal residual conditions on $\vx_m(t)$ and $\underline\vx_m(t)$ are equivalent to the minimal error conditions   
\[
\|\ve_m(t) \|_{A+tI}\to \min_{\vx_m(t)\in\spn(\vV_m)}\quad \text{and} \quad \|\underline\ve_m(t) \|_{A+tI}\to \min_{\underline\vx_m(t)\in\spn(V_m)}.
\]
Hence, 
\[
    \|\ve_m(t) \|_{A+tI} \leq \|\underline\ve_m(t) \|_{A+tI}
\]
as the minimization criterion defining the left-hand side takes place over the augmented subspace which contains the non-augmented Krylov space. 
Now, we are ready to conclude
\begin{eqnarray*}
    \| f(A)\vb - \vf_m \|_A &=& \left\| \int_0^{+\infty} \ve_m(t) \d\mu(t) \right\|_A \leq  \int_0^{+\infty} \|\ve_m(t)\|_A \d\mu(t) \\
    &\leq&  \int_0^{+\infty} {\frac{\sqrt{\lambda_{\max}}}{\sqrt{\lambda_{\max} + t}}} \cdot \|\ve_m(t)\|_{A + tI} \d\mu(t) \\
    &\leq&  \int_0^{+\infty} {\frac{\sqrt{\lambda_{\max}}}{\sqrt{\lambda_{\max} + t}}} \cdot \|\underline\ve_m(t)\|_{A + tI} \d\mu(t).
\end{eqnarray*}
The last expression is an upper bound involving the error $\underline\ve_m(t)$ of the non-augmented (standard) FOM approximant  $\underline\vx_m(t)$. It is precisely the expression which is further bounded in \cite{FrommerGuettelSchweitzer2014b} and the Lemma~4.1 and Corollary~4.4 therein. In summary, we obtain the following version of \cite[Cor.~4.4]{FrommerGuettelSchweitzer2014b}.

\begin{theorem}
Let $A\in\mathbb{C}^{N\times N}$ be Hermitian positive definite, $\vb\in\mathbb{C}^N$, and let $f(z) = \int_{t_0}^{+\infty} (z+t)^{-1}\d \mu(t)$ be a Stieltjes function ($t_0\geq 0)$. Let $\vf_m$ be the augmented FOM approximant as defined in \eqref{eq:afom2} and with $\vV_m = [ V_m, U]$,  where the columns of $V_m$ span $\mathcal{K}_m(A,\vb)$ and $U$ spans an augmentation space. 
Further, let $\lambda_{\min}$ and $\lambda_{\max}$ denote the smallest and largest eigenvalues of $A$, respectively, and define the functions
\[
\kappa(t) :=\frac{\lambda_{\max}+t}{\lambda_{\min}+t},
\quad
c(t) := \frac{\sqrt{\kappa(t)}-1}{\sqrt{\kappa(t)}+1}, 
\quad \text{and} \quad
\alpha_m(t) = \frac{1}{\cosh(m \ln c(t))}.
\]
Then the augmented FOM approximant satisfies
\[
\| f(A)\vb - \vf_m \|_A \leq C \alpha_m(t_0)
\]
with a constant $C = \|\vb\|_2\sqrt{\lambda_{\max}} \cdot f(\sqrt{\lambda_{\min}\lambda_{\max}})$. 
\end{theorem}

\smallskip

Note that this theorem guarantees exactly the same convergence as \cite[Cor.~4.4]{FrommerGuettelSchweitzer2014b} for the non-augmented (standard) FOM. Without adding further conditions on the augmentation space, which will likely be difficult to verify in practice, it appears to be unclear how to quantify any convergence improvements due to the augmentation. For an informal explanation of why convergence improvements can be expected when the augmentation space contains (approximate) eigenvectors of $A$ so that spectral deflation takes place, we refer to \cite[Sec.~4]{EiermannErnstGuettel2011}.


\subsection{Closed-form augmented GMRES-type approximation}

As mentioned previously, a GMRES-type approximant is obtained from \eqref{eq:fm} if $\vW_m = A \vV_m$. In the case where $\vV_m = [V_m]$ is an orthonormal basis of the Krylov space $\mathcal{K}_m(A,\vb)$ satisfying an Arnoldi decomposition \eqref{eq:arn}, \eqref{eq:fm} can be written as
\begin{eqnarray*}
\vf_m &=& V_m f\left( (V_m^* A^* V_m)^{-1}  V_m^* A^* A V_m   \right)  (V_m^* A^* V_m)^{-1} (A V_m)^* \vb\\
&=& V_m f\left(  H_m + (|h_{m+1,m}|^2 H_m^{-*}\ve_m) \ve_m^T \right) V_m^* \vb.
\end{eqnarray*}
This expression is well known (see, e.g., \cite[eq.~(1.5)]{GoossensRoose1999} and \cite[eq.~(2.5)]{HH05}) and it has been used  in \cite{FrommerGuettelSchweitzer2014b} 
(albeit with a typo) 
to show convergence of the (restarted) harmonic Arnoldi method for positive real matrices.

\medskip

\section{Randomized sketching}\label{sec:sketch}
As the dimension $m$ of the Krylov basis $V_m$ grows large, the cost of its orthogonalization may dominate. In order to overcome this drawback, we may work with a basis that is only partially orthogonalized and then use randomized sketching  to deal with the non-orthogonality. 
The representation \eqref{eq:afom2} is particularly well suited for sketching; see, e.g., \cite{martinsson2020randomized,balabanov2019randomized,balabanov2021randomized,BalabanovGrigori21,NakatsukasaTropp21,balabanov2022randomized}. 
Assume that we have a matrix $S\in\mathbb{C}^{s\times N}$ with $m < s\ll N$ which acts as an approximate isometry for the Euclidean norm $\|\,\cdot\,\|$. More precisely, given a positive integer $m$ and some $\varepsilon\in [0,1)$, let $S$ be such that for all vectors~$\vv\in\spn(\vV_m)$,
\begin{equation}
(1-\varepsilon) \| \vv \|^2 \leq \| S \vv\|^2 \leq (1+\varepsilon) \|\vv\|^2.
\label{eq:sketch}
\end{equation}
The mapping $S$ is called an \emph{$\varepsilon$-subspace embedding} for $\spn(\vV_m)$; see, e.g.,~\cite{sarlos2006improved,woodruff2014sketching,martinsson2020randomized}. 
Condition~\eqref{eq:sketch} can equivalently be stated with the Euclidean inner product~\cite[Cor.~4]{sarlos2006improved}: for all $\vu,\vv \in \spn(\vV_m)$,  
\begin{equation}{\label{eq:sketch_innerproduct}}
\langle \vu, \vv \rangle - \varepsilon \| \vu\|\cdot \|\vv\|
            \leq \langle S\vu, S\vv \rangle 
            \leq \langle \vu, \vv \rangle + \varepsilon \| \vu\|\cdot \|\vv\|.\nonumber
\end{equation}
 In practice, $S$ is not explicitly available but we can draw it at random to achieve~\eqref{eq:sketch} with high probability.

Using the sketching operator~$S$, we can replace all the implicit least squares problems in \eqref{eq:afom2} by solutions of their sketched counterpart, leading to 
\begin{equation}
\widehat \vf_m := \vV_m  f\left([S\vV_m]^\dagger [S A \vV_m]   \right)  [S\vV_m]^\dagger S\vb.
\label{eq:afom2sketch}
\end{equation}
This representation may be badly affected by ill-conditioning of the sketched basis~$S\vV_m$ and it is usually beneficial to perform basis whitening~\cite{rokhlin2008fast} by computing an economic QR decomposition $S \vV_m = Q_m R_m$ and using
\begin{equation}
\widehat \vf_m = \vV_m \left( R_m^{-1}  f\left(Q_m^* [S A \vV_m] R_m^{-1}   \right)  Q_m^*  S\vb\right).
\label{eq:afom2sketchwhite}
\end{equation}

It is interesting to note that \eqref{eq:afom2sketchwhite}  is formally the same as the sketched FOM approximant labeled (sFOM{'}{'}{'}) in \cite{GS23}. However, in that work $\vV_m$ was chosen as a basis of the Krylov space $\mathcal{K}_m(A,\vb)$, while here we allow for a basis of an arbitrary subspace of $\mathbb{C}^N$.  

\begin{remark}
The general form of the approximant \eqref{eq:afom} lends itself to sketching also in the case where $\vW_m = A\vV_m$, i.e., for the GMRES-type approximant. A~closed-form  sketched GMRES-type approximant is 
\begin{equation}
\widetilde{\vf}_m = \vV_m f\left([(S\vW_m)^* S\vV_m]^{-1} (S\vW_m)^* S\vW_m\right) [(S\vW_m)^* S\vV_m ]^{-1} (S\vW_m)^*S\vb.
\label{eq:sgmres}
\end{equation}
Note that all the sketched matrices appearing within $f(\cdot)$ are small (unrelated to the original problem dimension~$N$) and this formula is indeed relatively cheap to evaluate. It may hence be an attractive alternative to the quadrature-based 
sGMRES  approximant presented in \cite{GS23}. However, care must be taken with the numerical implementation as the condition number of the inverted matrices can become large. We have only performed some preliminary tests with the sketched GMRES-type approximant and empirically found it to work well, but here we prefer not to discuss this approximant any further and instead focus on sketched and recycled~FOM.
\end{remark}

\subsection{Updating the augmentation space}\label{sec:sRR}
We now turn our attention to the problem of computing a sequence of matrix functions $f(A^{(i)})\vb^{(i)}$ and discuss how to update the augmentation space $U$ from one problem to the next. Say, for one problem $f(A^{(i)})\vb^{(i)}$ we have computed a Krylov basis $V_m$ of $\mathcal{K}_m(A^{(i)},\vb^{(i)})$ and we have an augmentation basis $U\in\mathbb{C}^{N\times k}$ accumulated over previous problems. Neither $V_m$ nor $U$ are now assumed to be orthonormal. 

Consider the augmented basis matrix  $\vV_m = [ U, V_m]$ and assume that we have the sketches $S \vV_m$ and $S A^{(i)} \vV_m$ at our disposal. To obtain an updated augmentation space $U$, we can now follow the sRR approach presented in \cite[Sec.~6.3]{NakatsukasaTropp21} by solving the problem 
\begin{equation}
    \widehat M := \arg \min_{M\in\mathbb{C}^{(m+k)\times (m+k)}} \| S A^{(i)} \vV_m - S\vV_m M \|_F.
    \label{eq:srr}
\end{equation}
A solution (of minimum Frobenius norm) is $\widehat M = (S\vV_m)^\dagger (S A^{(i)} \vV_m)$ or, exploiting the basis whitening factorization $S \vV_m = Q_m R_m$ from above,  
\[
\widehat M = R_m^{-1} Q_m^* (S A^{(i)} \vV_m).
\]
We  can now compute a partial Schur decomposition 
\[
\widehat M X = X T, \quad X\in\mathbb{C}^{(m+k)\times k}, \quad T \in \mathbb{C}^{k\times k},
\]
associated with appropriately selected eigenvalues  of $\widehat M$. (If a quasi-triangular real Schur form is used, it may be necessary to increase $k$ to $k+1$ to not tear apart any $2\times 2$ diagonal blocks of $T$.) Finally, the augmentation space and its sketched version are updated as
\[
 U \leftarrow \vV_m X, \quad SU \leftarrow (S\vV_m) X.
\]
In addition, if the matrix $A^{(i)}=A^{(i+1)}$ remains unchanged for the next problem, we can also update $SA^{(i+1)} U \leftarrow (SA^{(i)}\vV_m) X$ without computing  any sketches or matrix-vector products with $A^{(i+1)}$.  
We summarize the overall FOM procedure which combines sketching and recycling in Algorithm~\ref{alg:srFOM}.

\begin{algorithm}[t]
\caption{Sketch-and-recycle FOM for  $f({A}^{(i)})\vb^{(i)}$ \hfill srFOM \label{alg:srFOM}}
\begin{algorithmic}[1]
\STATE{\textbf{Input:} ${A}^{(i)}$, $\vb^{(i)}$, integers $m,k$, sketch parameter $s$ ($m < s \ll N$)}
\STATE{\textbf{Output:} Approximants $\widehat{ \vf}_{m}^{(i)} \approx f({A}^{(i)}) \vb^{(i)}$}\\[1mm]
\STATE{Draw a sketching matrix $S \in \mathbb{C}^{s \times N}$}
\STATE{Initialize $U\leftarrow [\ ]$}
\FOR{$i = 1,2,\dots$}
\STATE{Generate (non-orthogonal) basis $V_{m}$ of $\mathcal{K}_{m}({A}^{(i)},\vb^{(i)})$} 
\STATE{Set $\vV_m \leftarrow [ V_m, U ]$} 
\STATE{Obtain sketch $S\vV_m$ (and $SA^{(i)}\vV_m$ if $i=1$ or $A^{(i)}\neq A^{(i-1)}$)}
\STATE{ Compute thin QR factorization ${Q}_{m} {R}_{m} \leftarrow {S} \vV_{m} $ }\\[1mm]
\STATE{$\widehat{\vf}_{m}^{(i)} \leftarrow \vV_{m} \left({R}_{m}^{-1} f({Q}^{*}_{m} [{S} {A}^{(i)} \vV_{m}] {R}_{m}^{-1}) {Q}^{*}_{m} {S} \vb^{(i)}\right)$}
\smallskip
\STATE{Compute $\widehat M \leftarrow R_m^{-1} Q_m^* [S A^{(i)} \vV_m]$}
\STATE{Compute partial Schur form $\widehat M X = X T$}
\STATE{Update  $U \leftarrow \vV_m X$, $SU \leftarrow S\vV_m X$ (and  $SA^{(i)}U \leftarrow SA^{(i)}\vV_m X$)}
\ENDFOR
\smallskip
\end{algorithmic}
\end{algorithm}

\section{Implementation}\label{sec:implem}

In this section we discuss a number of points relating to the practical implementation of Algorithm~\ref{alg:srFOM}.

\subsection{Stabilization}\label{sec:stabilized_sRR}
If the matrix $S \vV_{m} = S[V_m,U] $ is poorly conditioned, then we can employ the regularization proposed in \cite{NakatsukasaTropp21}. Instead of computing the  QR factorization $S \vV_m = Q_m R_m$, we  compute an economic SVD $ S \vV_{m} = L \Sigma J^{*}$. We define the truncated matrices $L_\ell=L(\,:\,,1:\ell)$, $\Sigma_\ell = \Sigma(1:\ell,1:\ell)$, $J_\ell= J(\,:\,,1:\ell)$, where $\ell$ is the largest integer such that the singular values of $S\vV_m$ satisfy
\[
\sigma_1\geq \sigma_2\geq \cdots \geq \sigma_\ell \geq \texttt{svdtol}
\]
for a given tolerance (like $\texttt{svdtol}=10^{-14}$). 
Replacing 
\[
    S \vV_{m} \rightarrow L_\ell, \quad 
    S A \vV_{m} \rightarrow S A \vV_{m}J_\ell \Sigma_\ell^{-1}, \quad 
    \vV_m \rightarrow \vV_m J_\ell \Sigma_\ell^{-1}   \ \ \text{(only implicitly!)}
\]
in \eqref{eq:afom2sketch}, we obtain the stabilized  sketched FOM approximant
\begin{equation}
\widehat \vf_m^{\text{stab}} = \vV_m \left( J_{\ell} \cdot \Sigma_{\ell}^{-1} f\big( L^{*}_{\ell} [S A \vV_{m}] J_{\ell} \Sigma^{-1}_{\ell}   \big)  L^{*}_{\ell}( S\vb)\right).
\label{eq:afom2sketchwhitestab}
\end{equation}
The recycling subspace can  be updated by computing an ordered QZ decomposition so that both of 
\[
   Q (L_\ell^{*} S A \vV_{m} J_\ell) Z,  \quad Q \Sigma_\ell Z
\]
are upper-triangular matrices with the $k$ targeted eigenvalues appearing in the upper-left block. (If a quasi-triangular real QZ form is used, it may be necessary to increase $k$ to $k+1$ to not tear apart any $2\times 2$ diagonal blocks.) The recycling subspace can then be updated as $U \leftarrow \vV_{m} J_\ell Z(\,:\, , 1\,: \, k)$.

\begin{remark}
The sketched GMRES-type approximant \eqref{eq:sgmres} could be stabilized in the same manner, but it is more convenient to compute a truncated SVD of $S\vW_m = SA\vV_m$ (instead of $S\vV_m$):  $S\vW_m \approx L_\ell \Sigma_\ell J_\ell^*$. Upon replacing 
\[
    S \vV_{m} \rightarrow S \vV_{m} J_\ell \Sigma_\ell^{-1}, \quad 
    S \vW_m \rightarrow L_\ell, \quad 
    \vV_m \rightarrow \vV_m J_\ell \Sigma_\ell^{-1}   \ \ \text{(only implicitly!)},
\]
the stabilized version of the approximant~\eqref{eq:sgmres} takes the form
\begin{equation*}
\widetilde{\vf}_m^{\text{stab}} = \vV_m \left(J_\ell   f\left( [L_\ell^* S\vV_mJ_\ell ]^{-1} \Sigma_\ell \right) 
 [L_\ell^* S\vV_mJ_\ell  ]^{-1} L_\ell^*(S\vb)\right).
\end{equation*}
\end{remark}

\subsection{Computational complexity}
We now  compare the cost of sketched-and-recycled FOM to sketched FOM for solving a single problem in the sequence \eqref{eq:sequence_of_matrix_function_applications}. Following the analysis in \cite{GS23} we assume that 
\smallskip
\begin{itemize}
    \item an Arnoldi cycle length of $m$ and a recycling subspace dimension of $k$ is used,
    \item $V_{m}$ is computed using the $t$-truncated Arnoldi process where $t = \mathcal{O}(1)$, and
    \item the sketching parameter is chosen as $ s = \mathcal{O}(m+k)$.
\end{itemize}
\smallskip
For non-Hermitian problems, the dominant arithmetic cost  of the Arnoldi method (aside the unavoidable matrix-vector products with $A$) are the $\mathcal{O}(N m^{2})$ arithmetic operations required to orthogonalize the  Krylov basis~$V_m$.  One of the main advantages of sFOM and srFOM is that these methods can work with non-orthogonal Krylov bases, e.g., generated by the $t$-truncated Arnoldi process. In that process, the $(j+1)$-st Krylov basis vector is computed by projecting $\vw_j := A \vv_{j+1}$ against the previous $t$ basis vectors
\begin{equation*}\label{eq:arntrunc}
\widehat \vw_j := \vw_j - \sum_{i=\max\{ 1, j+1-t\}}^j h_{i,j} \vv_i, \quad h_{i,j} := \vv_i^* \vw_j,
\end{equation*}
and then setting $\vv_{j+1}:=\widehat \vw_j/\|\widehat \vw_j\|$; see, e.g.,~\cite[Sec.~3.3]{Saad1981}.  
Truncated Arnoldi requires only $\mathcal{O}(N m t)$ arithmetic operations for the basis generation, i.e., the cost grows linearly in the basis dimension~$m$. In a recent work \cite{GS23b}, the possibility of using the sketched basis $SV_m$ (which is readily available) to select $t$ vectors among \emph{all} previously computed vectors (not just the most recent ones) has been explored. Other possibilities for constructing Krylov bases using a limited number of inner products (or no inner products at all) exist as well, e.g., based on recurrence relations for Chebyshev polynomials~\cite[Sec.~4]{JoubertCarey1992} or Newton polynomials~\cite[Sec.~4]{PhilippeReichel2012}. 

If the sketching operator $S$ is a  subsampled random discrete cosine or Fourier transform \cite{WoolfeLibertyRokhlinTygert2008,martinsson2020randomized}, then the cost of sketching a matrix $\vV_{m}$ with $m+k$ columns is  $\mathcal{O}(N(m+k) \mbox{log}(m+k))$. If $\vV_m = [U,V_m]$ and $V_m$ is constructed by the truncated Arnoldi process,  we can cheaply obtain the sketch  $S A V_{m}$ using the relation
\[
    S A V_{m} = (S V_m) H_m + h_{m+1,m} (S\vv_{m+1}) \ve_m^T,
\]
so overall only $m$ vector sketches are needed. (This observation also applies to sFOM~\cite{GS23} without recycling but has not been utilized in that work.) 

Performing the thin QR factorization in line $9$ requires a total cost of $ \mathcal{O}(s(m+k)^{2}) = \mathcal{O}((m+k)^{3})$, while for large enough $N$ forming the full approximant in line $10$ is dominated by the $O(N(m+k))$ cost for the linear combination with the columns of~$\vV_{m}$. This computation should hence be avoided whenever possible, working only with short $(m+k)$-dimensional  coefficient vectors $\left({R}_{m}^{-1} f({Q}^{*}_{m} [{S} {A} \vV_{m}] {R}_{m}^{-1}) {Q}^{*}_{m} {S} \vb\right)$ for tasks like error estimation (see the next subsection). The dominant cost of updating the recycling subspace (line~$13$) is in the $\mathcal{O}(N(m+k)k)$ arithmetic operations required to construct $U$. The full dominant cost of sketched-and-recycled FOM is thus $\mathcal{O}(N(m+k)\log(m+k) + (m+k)^{3})$. When compared to the $\mathcal{O}(Nm \log(m) + m^{3})$ arithmetic operations required in sketched FOM \cite{GS23}, we see that the dominant cost of both methods is the same when we choose an Arnoldi cycle length of $m' = m - k$.

\subsection{Error estimation}\label{sec:error_estimate}

It is widely appreciated that the convergence analysis of Krylov methods for nonsymmetric matrices is challenging. The difficulties will be even more pronounced when recycling and sketching techniques are incorporated into the algorithms. We do not currently have a practical error analysis of our sketched-and-recycled FOM (or GMRES) algorithms. Even in the case of sketched FOM or GMRES alone, existing error bounds are rarely useful in practice as they tend to overestimate the actual error by a large margin and often require information  about the matrix (like the numerical range) that is not easily available~\cite{GS23}. See also \cite{palitta2023sketched} for some interesting insights into sFOM.

Luckily, simple and practical a-posteriori error estimates are easy to derive, such as the frequently used difference of two iterates, i.e.,
\[
\big\|f(A)\vb-\widehat{\vf}_m\big\| \approx \big\|\widehat{\vf}_{m+d}-\widehat{\vf}_m\big\|
\]
for a small integer $d \geq 1$. Following \cite{GS23}, it is possible to approximately evaluate this stopping criterion without access to the full matrix $\vV_m$ and without forming  $\widehat{\vf}_{m+d}$ and $\widehat{\vf}_{m}$ explicitly. Using the assumption that  $S$ is an $\varepsilon$-subspace embedding for $\vV_m$, satisfying~\eqref{eq:sketch}, we have 
\begin{equation}\label{eq:sketched_err_est1}
{\frac{1}{\sqrt{1+\varepsilon}} \big\|S(\widehat{\vf}_{m+d}-\widehat{\vf}_m)\big\|}  \leq 
\big\|\widehat{\vf}_{m+d}-\widehat{\vf}_m\big\| \leq \frac{1}{\sqrt{1-\varepsilon}}\big\|S(\widehat{\vf}_{m+d}-\widehat{\vf}_m)\big\|
\end{equation}
Writing $\widehat{\vf}_{m} = \vV_m\widehat{\vy}_m$, we obtain from~\eqref{eq:sketched_err_est1} the bound
\begin{eqnarray*}
\big\|\widehat{\vf}_{m+d}-\widehat{\vf}_m\big\| &\leq& \frac{1}{\sqrt{1-\varepsilon}}\big\|SV_{m+d}\widehat{\vy}_{m+d}-SV_m\widehat{\vy}_m\big\|\\ 
&=& \frac{1}{\sqrt{1-\varepsilon}}\left\| SV_{m+d}\left(\widehat{\vy}_{m+d}-\left[\begin{array}{c} \widehat{\vy}_{m} \\ \vnull_d \end{array}\right]\right)\right\|.
\end{eqnarray*}
For an estimator, the unknown embedding quantity $\varepsilon$  can be  set to a  fixed generous constant (like $\varepsilon=0.99$), or it can be estimated by keeping track of $\|S\vv_j\|/\|\vv_j\|$ as the Krylov basis vectors $\vv_j$ are generated ($j = 1,\dots,m$).

\section{Numerical experiments}\label{sec:numex}

In this section we present results of numerical experiments demonstrating the effectiveness of the rFOM (Algorithm~\ref{alg:rFOM}) and srFOM (Algorithm~\ref{alg:srFOM}) introduced here. We compare these methods to the standard Arnoldi approximation (FOM), the quadrature-based rFOM (quad)~\cite[Algorithm~2]{burke2022krylov}, and sketched FOM (sFOM)~\cite[Algorithm~1]{GS23}. 

In all methods where recycling is used, the recycling subspace is updated with $k$~approximate eigenvectors corresponding to  eigenvalues closest to the origin.  {In all cases where randomized sketching is employed, the sketching operator is a subsampled randomized discrete cosine transform (see, e.g.,~\cite[Sec.~9.3]{martinsson2020randomized}) and the Arnoldi truncation parameter is $t=2$.} 
The sFOM algorithm uses basis whitening, but we run it \emph{without any stabilization} to provoke potential instabilities. Stabilization \emph{is}  used in srFOM (stab) with a fixed tolerance of $\texttt{svdtol}= 10^{-14}$. 

All experiments were performed in MATLAB\footnote{ Code available at \url{https://github.com/burkel8/srFOM}} R2023a on a Windows~11 HP laptop with an 11th Gen Intel(R) Core processor with 2.80\,GHz and 8\,GB of RAM. All reported runtimes are averages over 10~repetitions of an experiment. 

\subsection{Inverse square root}\label{sec:qcd}
In this set of experiments we consider the  function $f(z) = z^{-1/2}$. This function arises in several applications, including Dirichlet-to-Neumann maps (see, e.g., \cite{DGK16a}) and quantum chromodynamics (QCD)~(see, e.g.,~\cite{knechtli2017lattice, burke2022krylov}). Numerical simulations of QCD and the computation of spectral projectors often require the evaluation of~\eqref{eq:second_solution_approximation} involving the matrix sign function,  which can be represented as $\sign(A)=f(A^2)A\vb$. 

We first approximate a sequence of $20$ vectors of the form \eqref{eq:sequence_of_matrix_function_applications}  where the first matrix $A^{(1)}$ in the sequence is a complex non-Hermitian  lattice QCD matrix of size $N=3,072$ obtained from the SuiteSparse Matrix Collection \cite{Davis2011} (matrix ID~$1597$), plus a shift of $6.0777\cdot I$ to move all eigenvalues away from the negative real axis (the left-most real  eigenvalue of the unshifted matrix is $-6.0776$). The other matrices in the sequence take the form
\[
    A^{(i+1)} = A^{(i)} + 10^{-8} M^{(i)}, \quad  i = 1,2, \ldots, 19,
\]
where each $M^{(i)}$ is a unit Gaussian random matrix with  the same sparsity pattern as~$A^{(1)}$, and  the vectors $\vb^{(i)}$ are randomly generated with unit Gaussian entries.

Let us first discuss Figure~\ref{fig:compare_quad_v_closed} from the introduction in some more detail. This figure shows the relative error 
\[
\texttt{relerr} = \frac{\|f(A^{(i)})\vb^{(i)} - \vf_m^{(i)}\|}{\|f(A^{(i)})\vb^{(i)}\|}
\]
obtained after $m = 100$ Arnoldi iterations of rFOM (\cref{alg:rFOM}) and srFOM (\cref{alg:srFOM}). A recycling subspace dimension of $k = 30$, sketching dimension of~$s = 400$, and an Arnoldi truncation parameter $t = 2$ are used.  We also show the  relative error obtained from the standard FOM approximation and two runs of rFOM (quad) \cite{burke2022krylov} with $50$ and $100$ quadrature nodes, respectively.
The error curves demonstrate the effectiveness of rFOM and srFOM in reducing the relative error as the sequence of problems progresses and the recycling subspaces improves. 
Additionally, Figure~\ref{fig:compare_quad_v_closed} demonstrates that rFOM is more numerically robust than rFOM (quad), particularly for larger values of~$m$. We believe this improved robustness is at least partially attributable to the explicit orthogonalization of~$\vV_{m}$. Thus for the remainder of these experiments, we no longer include rFOM (quad) in the comparisons.

\begin{figure}
    \centering
    \includegraphics[width=.8\textwidth]{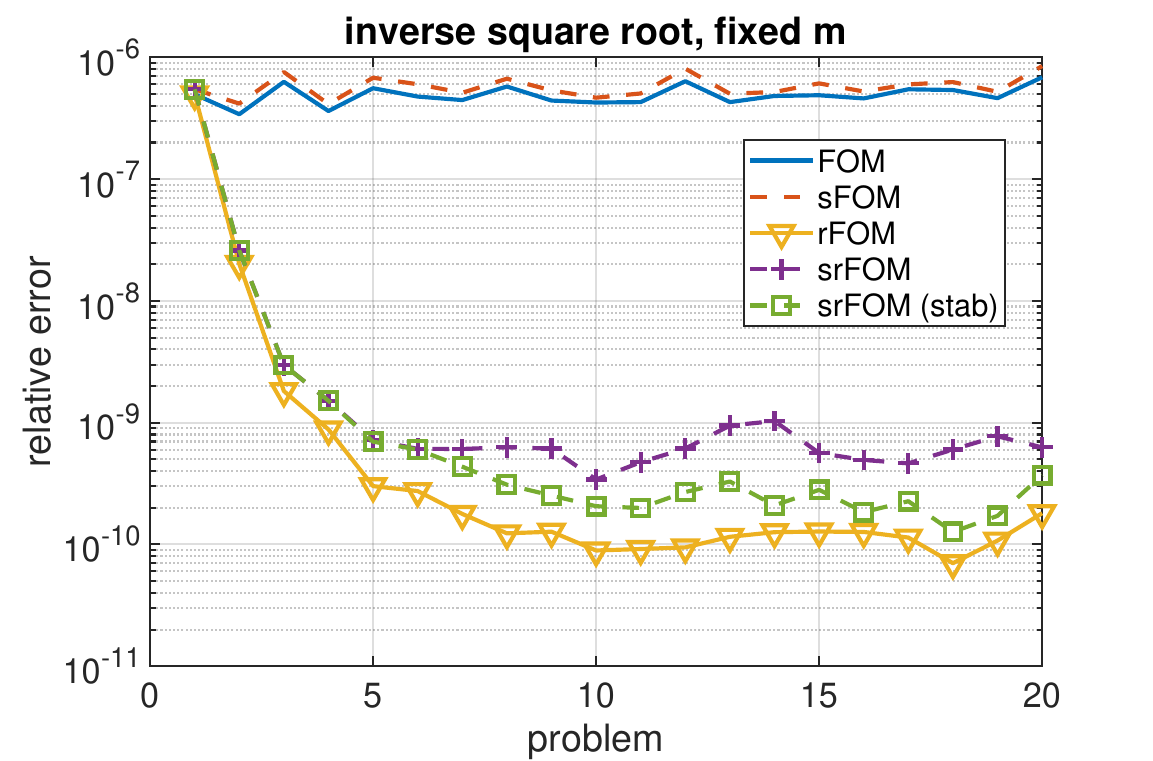}
    \caption{Relative errors when approximating  a sequence of $20$ vectors of the form $(A^{(i)})^{-1/2} \vb^{(i)}$ using FOM, sFOM, rFOM, srFOM, and stabilized srFOM. The Arnoldi cycle length is fixed at $m=100$ for all 20~problems. The other parameters are $k = 30$ (recycling subspace dimension), $s = 400$ (sketching dimension), and $t = 2$ (Arnoldi truncation length).}
    \label{fig:inv_sqrt_fixed_m.pdf}
\end{figure}

We now repeat the same experiment for Figure~\ref{fig:inv_sqrt_fixed_m.pdf}, but also include the  error curves for sFOM and the stabilized srFOM. We see that the sFOM approximants attain an error that is only slightly larger than that of FOM, and that stabilization indeed helps to further reduce the error of srFOM. 

\begin{figure}
    \centering
    \includegraphics[width=.8\textwidth]{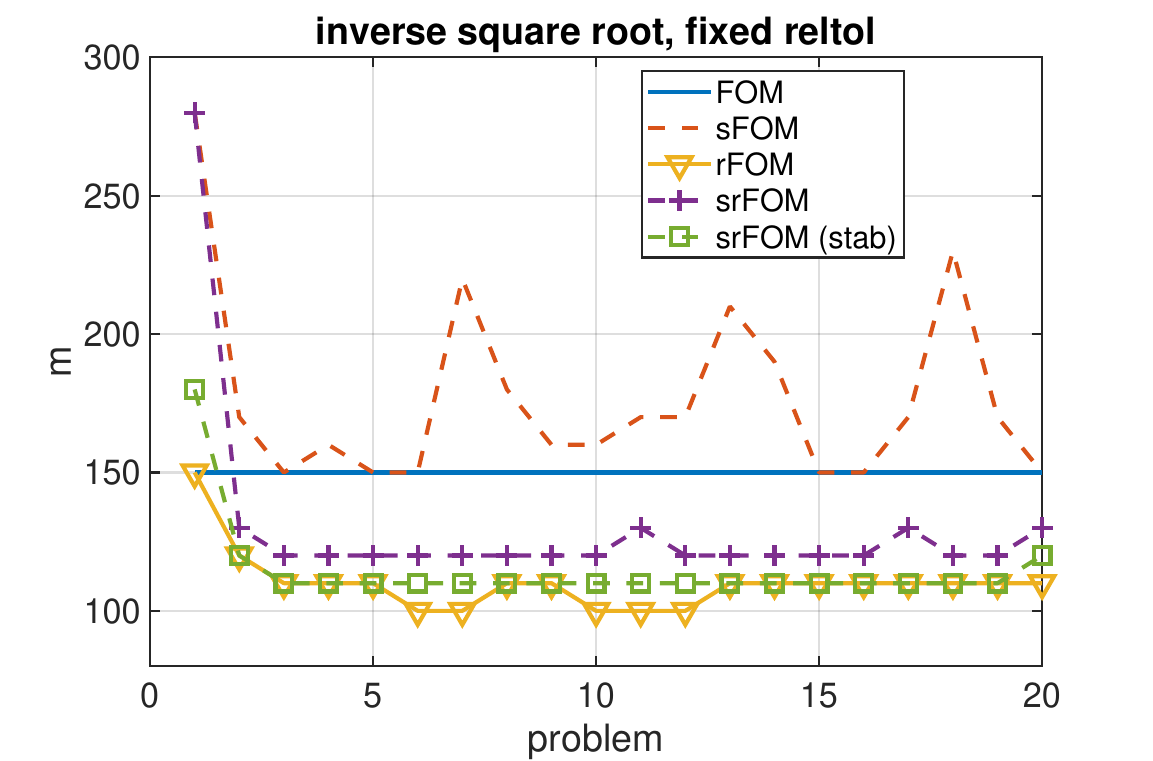}
    \caption{
    Arnoldi cycle length $m$ needed to approximate a sequence of $20$ vectors $(A^{(i)})^{-1/2} \vb^{(i)}$. Here the Krylov dimension~$m$ is chosen adaptively for each problem to reach a relative  error below $\texttt{reltol} = 10^{-10}$. The error criterion is checked every $d=10$ iterations. The other parameters are $k = 30$, $s = 400$, and $t = 2$. The rather irregular behavior of sFOM (and srFOM on the first problem) is caused by the lack of stabilization.
    }
    \label{fig:inv_sqrt_adaptive_m}
\end{figure}

In Figure~\ref{fig:inv_sqrt_adaptive_m} we now allow the number of Arnoldi iterations $m$ to vary for each problem in the sequence. This is done by checking the relative error every $d=10$ iterations and stopping when $\texttt{relerr}$ is below $10^{-10}$. We then plot the values of $m$ required for each problem. We first note the rather high value of~$m$ required by srFOM for the first problem and  the irregular behavior of sFOM for many of the problems. This is caused by the absence of  truncated SVD-based stabilization. (On the first problem where there is no available augmentation space~$U$, both sFOM and srFOM are equivalent.) The rather stringent tolerance of $10^{-10}$ pushes the non-stabilized sketching methods to the limit.  If the tolerance is raised to $10^{-9}$ or stabilization is employed, the irregular behavior goes away. This can be seen in Figure~\ref{fig:inv_sqrt_adaptive_m} for the srFOM (stab) method, which behaves much more regularly with $m$ steadily decreasing as the sequence of problems progresses. Indeed, the values of $m$ required by srFOM (stab) are only marginally above those of rFOM (which uses no sketching).

In Table~\ref{tab:QCD_stats} we record the total number of matrix-vector products with the matrices~$A^{(i)}$ (MAT-VEC's), the number of inner products, the number of vector sketches, and the total  runtime. 
We see that sFOM and rFOM are able to reduce the number of inner products and MAT-VEC's when compared to FOM, respectively, resulting in a moderate reduction of total runtime. However, a much better overall performance is obtained with the \emph{combination} of recycling and sketching (with stabilization). This is because srFOM (stab) combines the best of both worlds: it reduces the number of MAT-VEC's using recycling and it reduces orthogonalization cost via sketching.

\begin{table}[t]
\centering
\begin{tabular}{|p{2.9cm}||p{1.3cm}|p{1.3cm} |p{1.3cm}|p{1.3cm}|p{2.2cm}|}
 \hline
 \multicolumn{6}{|c|}{inverse square root, fixed \texttt{reltol}} \\
 \hline \hline
   & FOM & sFOM & rFOM & srFOM &  srFOM (stab)\\
 \hline \hline
    MAT-VEC's  & 3,000 & 3,610 & 2,660 & 3,140 & 2,860 \\
    Inner Products  & 229,500 & 10,810 & 942,285 & 7,690  & 6,850  \\
    Sketches  & 0 & 3,630 & 0 & 3,160 & 2,880 \\
    Runtime (seconds) &   10.6  &  9.8 & 8.4  &   6.6 &  5.7 \\ 
    \hline 
\end{tabular}
 \caption{
 Computational resources  required to approximate a sequence of $20$ vectors $(A^{(i)})^{-1/2} \vb^{(i)}$ to reach a relative error below  $\texttt{reltol} = 10^{-10}$. The error criterion is checked every $d=10$ iterations. The other parameters are $k = 30$, $s = 400$, and $t = 2$.
  }
 \label{tab:QCD_stats}
\end{table}

\subsection{Linear systems of equations}
Many applications in scientific computing require the solution to a sequence of slowly changing linear systems of the form 
\[
A^{(i)} \vx^{(i)} = \vb^{(i)},\qquad    i = 1,2,\ldots,
\]
a special case of \eqref{eq:sequence_of_matrix_function_applications} with $f(z) = z^{-1}$. Such problems arise, for example, with PDE-constraint optimization problems \cite{eveolving_structures}.  

Although Krylov subspace recycling is a well established technique for solving sequences of linear systems (see, e.g., \cite{parks2006recycling, eveolving_structures, gaul2014recycling}),  srFOM  appears to be the first attempt to combine recycling with randomized sketching. Of course, an efficient Krylov solver for linear systems should also incorporate some preconditioning and potentially restarting techniques. We will leave this to future work but still think it is worth demonstrating that srFOM is a good starting point for developing efficient linear system solvers that combine recycling and  sketching.

In this experiment we solve a sequence of $30$ linear systems with randomly generated right-hand sides $\vb^{(i)}$ having unit Gaussian entries. The matrix~$A$ is the \texttt{Neumann} matrix in MATLAB's matrix gallery of size $N=10,609$, which we shift by $+0.001\cdot I$ to make it nonsingular. 

\begin{figure}
    \centering
    \includegraphics[width=.8\textwidth]{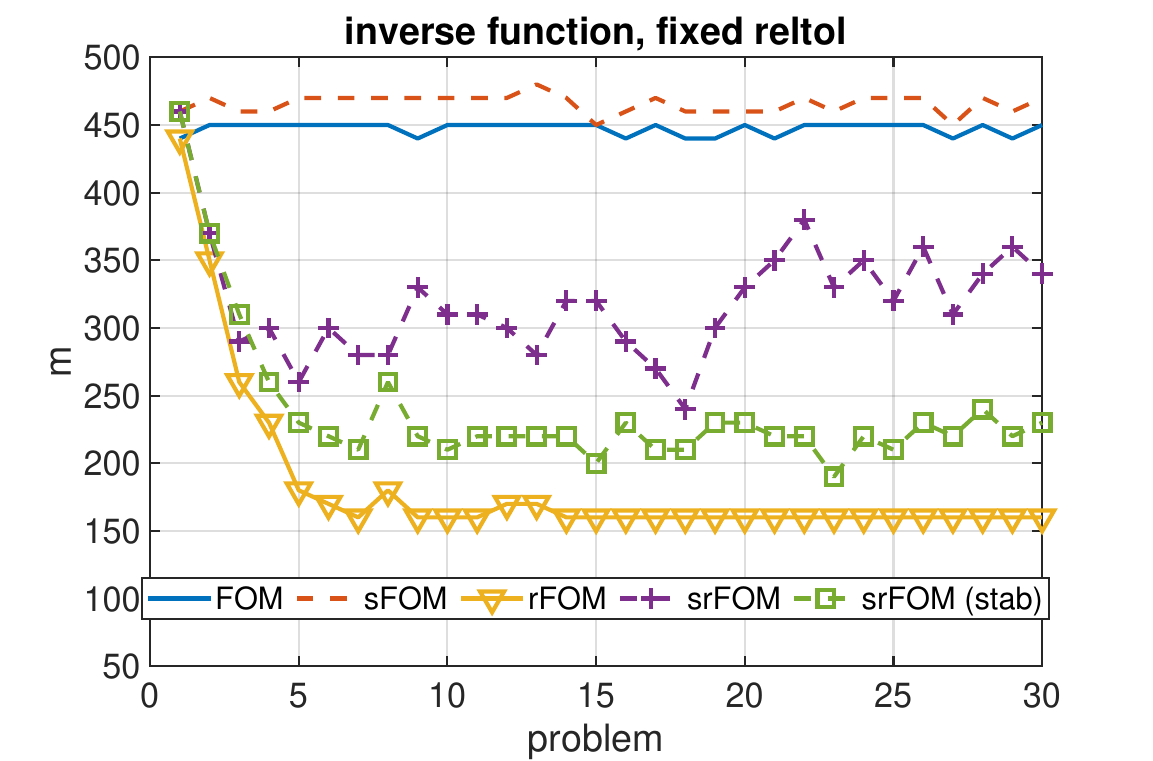}
    \caption{
     Arnoldi cycle length $m$ needed to approximate a sequence of $30$ vectors of the form $(A^{(i)})^{-1} \vb^{(i)}, i=1,2,\ldots,30$. Here the Krylov dimension~$m$ is chosen adaptively for each problem to reach a relative  error below $\texttt{reltol} = 10^{-9}$. The error criterion is checked every $d=10$ iterations. The other parameters are $k = 30$, $s = 900$, and $t = 2$.}
    \label{fig:inv_adaptive_m}
\end{figure}

In Figure~\ref{fig:inv_adaptive_m} we plot the Arnoldi cycle length~$m$ required for FOM, sFOM, rFOM, srFOM and srFOM (stab) to solve each linear system in the sequence with a relative error below $10^{-9}$. (To be consistent with the rest of the paper, we use an error-based stopping criterion, not a residual-based one.) The stopping criterion is monitored every $d = 10$ iterations. All recycling methods use a recycling subspace dimension of size $k = 30$, and all sketching methods use a sketching dimension of~$s = 900$ and an Arnoldi truncation parameter~$t = 2$. It is clear from Figure~\ref{fig:inv_adaptive_m} that rFOM and srFOM require significantly fewer Arnoldi iterations to reach convergence than FOM and sFOM. Additionally, Figure~\ref{fig:inv_adaptive_m} again illustrates the importance of stabilization, with srFOM (stab) requiring  consistently fewer Arnoldi iterations than srFOM. 

In Table~\ref{tab:Neumann_stats} we record the total number of matrix-vector products with the matrix~$A$ (MAT-VEC's), the number of inner products, the number of vector sketches,  and the total  runtime to solve the sequence of 30~problems. Note that rFOM alone already leads to a significant reduction of MAT-VEC's and  runtime over FOM, though the number of required inner products for the orthogonalization is very large. This is mitigated when sketching is employed, with srFOM (stab) again being the fasted method overall. 

\begin{table}
\centering
\begin{tabular}{|p{2.9cm}||p{1.3cm}|p{1.3cm} |p{1.3cm}|p{1.3cm}|p{2.2cm}|}
 \hline
 \multicolumn{6}{|c|}{inverse function, fixed \texttt{reltol}} \\
 \hline \hline
   & FOM & sFOM & rFOM & srFOM &  srFOM (stab)\\
 \hline \hline
    MAT-VEC's   & 13,420 & 13,970 &  5,510 & 9,580  & 7,140 \\
    Inner Products   & 3,022,030 & 41,880 & 3,975,815  & 28,710 & 21,390 \\
    Sketches  & 0 & 14,000 & 0 & 9,610 & 7,170 \\
    Runtime (seconds) &  84.2  &  42.7 & 29.8 & 25.8  &  22.0  \\ 
    \hline 
\end{tabular}
 \caption{
 Computational resources  required to approximate a sequence of $30$ vectors $(A^{(i)})^{-1} \vb^{(i)}$ to reach a relative error below  $\texttt{reltol} = 10^{-9}$. The error criterion is checked every $d=10$ iterations. The other parameters are $k = 30$, $s = 900$, and $t = 2$.
  }
 \label{tab:Neumann_stats}
\end{table}

\subsection{Exponential function}

We now consider the function  $f(z) = \exp(z)$ which arises, e.g., with the exponential time integration of ordinary differential equations. While Krylov recycling is most suited for functions with finite singularities, we demonstrate that the methods introduced in this paper can be readily applied to entire functions as well. We include this example also because the exponential has been a particularly difficult function to work with in quadrature-based restarting and recycling methods as there is no canonical contour for the integral representation~\cite{FrommerGuettelSchweitzer2014a,burke2022krylov}. 

For our test we use the FEM discretization of an advection--diffusion problem described in~\cite[Section~6.3]{AfanasjewEtAl2008a}. The matrix~$A$ of size $N=2,157$ was generated using the {COMSOL Multiphysics}
 software. We consider a sequence of 30~problems $\exp(0.01\cdot A)\vb^{(i)}$ where $\vb^{(1)}$ is chosen at random with unit Gaussian entries and every subsequent vector in the sequence is chosen as $\vb^{(i+1)} = \exp(0.01\cdot A)\vb^{(i)}$. This mimics an exponential time stepping method for the solution of $\vu'(t) = A \vu$, $\vu(0) = \vb^{(1)}$. The dimension of the augmentation space is $k=50$, the sketching parameter is $s=400$, and the Arnoldi truncation parameter is $t=2$.

\begin{figure}
    \centering
\begin{minipage}{.49\textwidth}
\hspace*{-3mm}\includegraphics[width=1.1\textwidth]{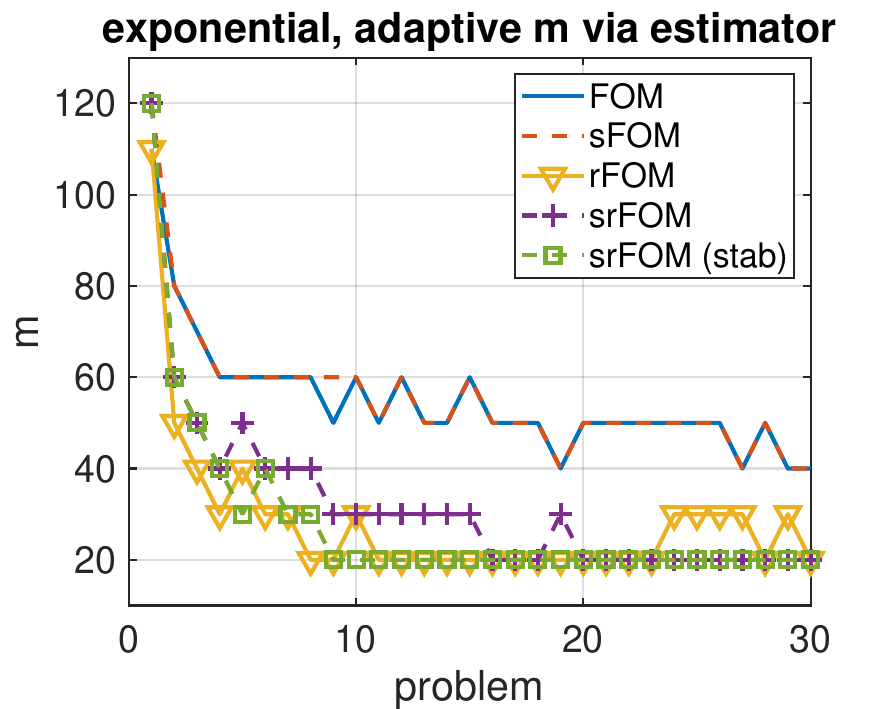}
\end{minipage}
\begin{minipage}{.49\textwidth}
\includegraphics[width=1.1\textwidth]{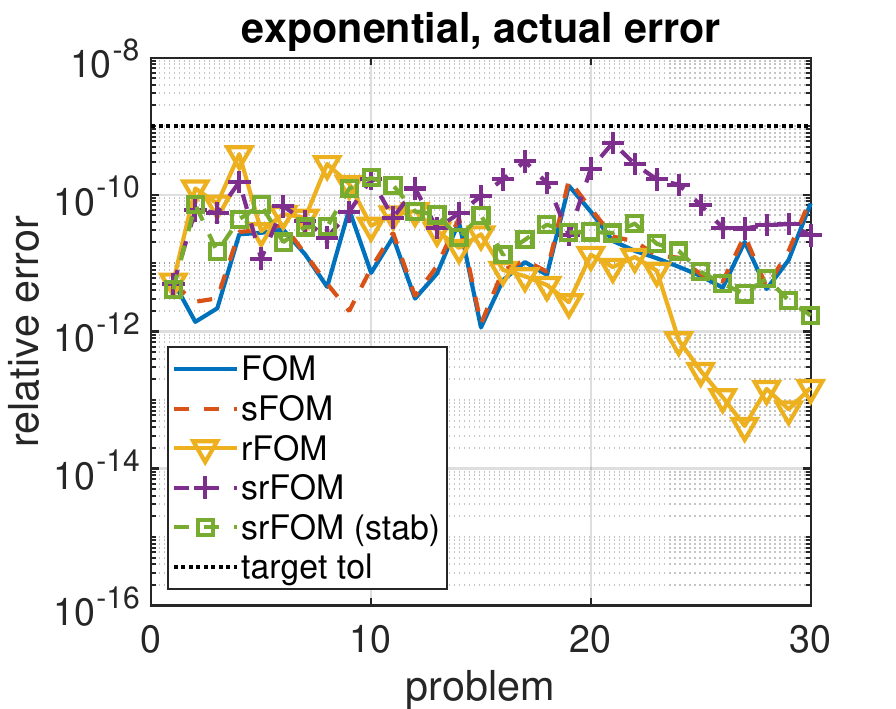}
\end{minipage}
    \caption{Arnoldi cycle length~$m$ (left) and relative errors for  approximating  a sequence of $30$~vectors  $\exp(0.01\cdot A)\vb^{(i)}$ with an adaptive control of $m$ using an error estimator with delay parameter $d=10$ and a target tolerance of $10^{-9}$. The other parameters are $k = 50$, $s = 400$, and $t = 2$.}
    \label{fig:exp_error_curves}
\end{figure}

In this test we use the error estimator described in section~\ref{sec:error_estimate} to terminate the Arnoldi iteration for each problem, checking the estimator every $d=10$ iterations for a relative error tolerance of $10^{-9}$. Figure~\ref{fig:exp_error_curves} (left) shows the number of Arnoldi iterations~$m$ performed for each problem. Note that $m$ decreases for all methods as the problems progress, even for the non-recycled FOM and sFOM. This is because the solution vectors $f(A)\vb^{(i)}$ get progressively easier to approximate for later time steps as certain eigenvector components in the starting vector~$\vb^{(i)}$ get damped out by the exponential $f(z)=\exp(0.01 z)$. Nevertheless, the recycling methods still achieve  a considerable further reduction in $m$. 
On the right of Figure~\ref{fig:exp_error_curves} we show the actual relative error of the computed approximants (we have used MATLAB's \texttt{expm} to compute the exact solutions, which is still feasible for this problem size). We find that these errors are indeed below the target tolerance of~$10^{-9}$ as intended, even though the error estimator may lead to a larger~$m$ than necessary in some cases. 

\section{Conclusions and future work}

We have presented a new Krylov subspace recycling method  for matrix functions (rFOM) which, in contrast to existing methods, does not require any numerical quadrature and performs more numerically robust. With srFOM we have shown how randomized sketching can easily be incorporated to further reduce the orthogonalization costs. 
We have presented results of numerical experiments which suggest that rFOM and srFOM can approximate a sequence of matrix functions with improved accuracy over FOM or sFOM, requiring less computational work and runtime.

One of the observations we have made but not discussed in the numerical experiments with recycling is that in some cases it is not actually necessary to form $A^{(i+1)} U$ nor $SA^{(i+1)} U$ explicitly without any delay in convergence even when $A^{(i+1)}\neq A^{(i)}$. In such cases the number of matrix-vector products per problem can be further reduced by~$k$, the dimension of the augmentation space. With such a modification, the sRR step~\eqref{eq:srr} for  problem $i+1$ in the sequence becomes  
\[
    \widetilde M := \arg \min_{M\in\mathbb{C}^{(m+k)\times (m+k)}} \| S  [A^{(i)} U^{(i)}, A^{(i+1)} V_m^{(i+1)}]  - S [U^{(i)}, V_m^{(i+1)}] M \|_F,
\]
where $U^{(i)}$ is the augmentation space extracted for the $i$-th problem and $V_m^{(i+1)}$ is a Krylov basis for $\mathcal{K}_m(A^{(i+1)},\vb^{(i+1)})$. This is an inexact Rayleigh--Ritz procedure for the matrix $A^{(i+1)}$ if $A^{(i)}\neq A^{(i+1)}$, but as the sketching with~$S$ introduces inexactness anyway, one might get away with it. The analysis of this is potential future work.

In future work we would like to explore how mixed-precision computations could be used to further reduce the runtime of rFOM and srFOM. Mixed-precision can be combined with recycling in multiple ways. For example, one might carry out parts of the Arnoldi process for forming $V_m$ using different precisions as in  \cite{BalabanovGrigori21}, one could only do the sketching in lower precision (similarly to \cite{georgiou2023mixed} which use a lower-precision sketch to obtain a preconditioner for LSQR), or one could use mixed-precision iterative refinement in the context of linear systems as in~\cite{oktay2022mixed}. Also, we believe that a thorough comparison of the proposed GMRES-type methods (with and without sketching) against more established recycling methods like GCRO-DR~\cite{parks2006recycling} for sequences (shifted) linear systems would be very interesting. 

For the case of linear systems of equations, we believe that srFOM and srGMRES can be easily combined with preconditioning. However, if the preconditioner is already effective in reducing the number of required MAT-VEC's and hence orthogonalization cost significantly, then the advantage of  randomized sketching dimishes. Hence sketching might be particularly useful in situations where an efficient preconditioner is not easily available; see also the discussion in~\cite[Sec.~1.4]{NakatsukasaTropp21}. We hope that such  aspects will be explored in future work.

\section*{Acknowledgments}
This work was jointly funded by an Irish Research Council Government Of Ireland Postgraduate Scholarship and the Manchester Mathematical Sciences (MiMS). We are grateful for discussions with Marcel Schweitzer and Kirk Soodhalter.

\bibliographystyle{plain}
\bibliography{refs}

\newcommand{\noopsort}[1]{} \newcommand{\printfirst}[2]{#1}
  \newcommand{\singleletter}[1]{#1} \newcommand{\switchargs}[2]{#2#1}
\begin{thebibliography}{10}

\bibitem{AfanasjewEtAl2008a}
M.~Afanasjew, M.~Eiermann, O.~G. Ernst, and S.~G{\"u}ttel.
\newblock Implementation of a restarted {K}rylov subspace method for the
  evaluation of matrix functions.
\newblock {\em Linear Algebra Appl.}, 429:2293--2314, 2008.

\bibitem{BalabanovGrigori21}
O.~Balabanov and L.~Grigori.
\newblock Randomized block {G}ram-{S}chmidt process for solution of linear
  systems and eigenvalue problems.
\newblock Technical Report arXiv:2111.14641, 2021.

\bibitem{balabanov2022randomized}
O.~Balabanov and L.~Grigori.
\newblock {Randomized Gram--Schmidt process with application to GMRES}.
\newblock {\em SIAM J. Sci. Comput.}, 44(3):A1450--A1474, 2022.

\bibitem{balabanov2019randomized}
O.~Balabanov and A.~Nouy.
\newblock Randomized linear algebra for model reduction. {Part I: Galerkin
  methods and error estimation}.
\newblock {\em Adv.\ Comput.\ Math.}, 45:2969--3019, 2019.

\bibitem{balabanov2021randomized}
O.~Balabanov and A.~Nouy.
\newblock Randomized linear algebra for model reduction--part {II}: minimal
  residual methods and dictionary-based approximation.
\newblock {\em Adv.\ Comput.\ Math.}, 47:1--54, 2021.

\bibitem{eveolving_structures}
M.~Bolten, E.~de~Sturler, C.~Hahn, and M.~L. Parks.
\newblock Krylov subspace recycling for evolving structures.
\newblock {\em Comput. Methods Appl. Mech. Engrg.}, 391:Paper No. 114222, 15,
  2022.

\bibitem{burke2022krylov}
L.~Burke, A.~Frommer, G.~Ramirez-Hidalgo, and K.~M. Soodhalter.
\newblock Krylov subspace recycling for matrix functions.
\newblock {\em arXiv preprint arXiv:2209.14163}, 2022.

\bibitem{exp_integrators}
M.~Caliari, F.~Cassini, and F.~Zivcovich.
\newblock B{AMPHI}: matrix-free and transpose-free action of linear
  combinations of {$\varphi$}-functions from exponential integrators.
\newblock {\em J. Comput. Appl. Math.}, 423:114973, 2023.

\bibitem{Davis2011}
T.~A. Davis and Y.~Hu.
\newblock The {U}niversity of {F}lorida sparse matrix collection.
\newblock {\em ACM TOMS}, 38:1--25, 2011.

\bibitem{DGK16a}
V.~Druskin, S.~G\"{u}ttel, and L.~Knizhnerman.
\newblock Near-optimal perfectly matched layers for indefinite {H}elmholtz
  problems.
\newblock {\em SIAM Review}, 58(1):90--116, 2016.

\bibitem{EiermannErnstGuettel2011}
M.~Eiermann, O.~G. Ernst, and S.~G\"uttel.
\newblock Deflated restarting for matrix functions.
\newblock {\em SIAM J.\ Matrix Anal.\ Appl.}, 32(2):621--641, 2011.

\bibitem{FrommerGuettelSchweitzer2014b}
A.~Frommer, S.~G\"uttel, and M.~Schweitzer.
\newblock Convergence of restarted {K}rylov subspace methods for {S}tieltjes
  functions of matrices.
\newblock {\em SIAM J.\ Matrix Anal.\ Appl.}, 35(4):1602--1624, 2014.

\bibitem{FrommerGuettelSchweitzer2014a}
A.~Frommer, S.~G{\"u}ttel, and M.~Schweitzer.
\newblock Efficient and stable {A}rnoldi restarts for matrix functions based on
  quadrature.
\newblock {\em SIAM J.\ Matrix Anal.\ Appl.}, 35:661--683, 2014.

\bibitem{frommer2000numerical}
A.~Frommer, T.~Lippert, B.~Medeke, and K.~Schilling.
\newblock {\em {N}umerical Challenges in {L}attice {Q}uantum {C}hromodynamics:
  {J}oint Interdisciplinary Workshop of John Von Neumann Institute for
  Computing, J{\"u}lich, and Institute of Applied Computer Science, Wuppertal
  University, August 1999}, volume~15.
\newblock Springer Science \& Business Media, 2000.

\bibitem{frommer2017block}
A.~Frommer, K.~Lund, and D.~B. Szyld.
\newblock {B}lock {K}rylov subspace methods for functions of matrices.
\newblock {\em {E}lectron. {T}rans. {N}umer. {A}nal.}, 47:100--126, 2017.

\bibitem{frommer2020block}
A.~Frommer, K.~Lund, and D.~B. Szyld.
\newblock {B}lock {K}rylov subspace methods for functions of matrices {II}:
  {M}odified block {FOM}.
\newblock {\em {SIAM} {J}. {M}atrix {A}nal. {A}ppl.}, 41(2):804--837, 2020.

\bibitem{gaul2014recycling}
A.~Gaul.
\newblock {\em {R}ecycling {K}rylov subspace methods for sequences of linear
  systems}.
\newblock PhD thesis, Technische Universität Berlin, Fakultät II - Mathematik
  und Naturwissenschaften, 2014.

\bibitem{georgiou2023mixed}
V.~Georgiou, C.~Boutsikas, P.~Drineas, and H.~Anzt.
\newblock A mixed precision randomized preconditioner for the {LSQR} solver on
  {GPU}s.
\newblock In {\em International Conference on High Performance Computing},
  pages 164--181. Springer, 2023.

\bibitem{GoossensRoose1999}
S.~Goossens and D.~Roose.
\newblock Ritz and harmonic {R}itz values and the convergence of {FOM} and
  {GMRES}.
\newblock {\em Numer.\ Linear Algebra Appl.}, 6(4):281--293, 1999.

\bibitem{GS23}
S.~G\"{u}ttel and M.~Schweitzer.
\newblock Randomized sketching for {K}rylov approximations of large-scale
  matrix functions.
\newblock {\em SIAM J. Matrix Anal. Appl.}, 44(3):1073--1095, 2023.

\bibitem{GS23b}
S.~G\"{u}ttel and I.~Simunec.
\newblock {A} sketch-and-select {A}rnoldi process.
\newblock Technical Report arXiv:2306.03592, 2023.

\bibitem{HH05}
M.~Hochbruck and M.~E. Hochstenbach.
\newblock Subspace extraction for matrix functions.
\newblock Technical report, Case Western Reserve University, Department of
  Mathematics, Cleveland, 2005.

\bibitem{JoubertCarey1992}
W.~D. Joubert and G.~F. Carey.
\newblock Parallelizable restarted iterative methods for nonsymmetric linear
  systems. part {I}: {T}heory.
\newblock {\em Int.\ J.\ Comput.\ Math.}, 44(1-4):243--267, 1992.

\bibitem{knechtli2017lattice}
F.~Knechtli, M.~G{\"u}nther, and M.~Peardon.
\newblock {\em {L}attice {Q}uantum {C}hromodynamics: {P}ractical {E}ssentials}.
\newblock Springer, 2017.

\bibitem{Lu2003}
Y.~Y. Lu.
\newblock {Computing a matrix function for exponential integrators}.
\newblock {\em J. Comput. Appl. Math.}, 161(1):203--216, 2003.

\bibitem{martinsson2020randomized}
P.-G. Martinsson and J.~A. Tropp.
\newblock Randomized numerical linear algebra: {F}oundations and algorithms.
\newblock {\em Acta Numer.}, 29:403--572, 2020.

\bibitem{Morgan2002}
R.~B. Morgan.
\newblock {GMRES} with deflated restarting.
\newblock {\em SIAM J.\ Sci.\ Comput.}, 24(1):20--37, 2002.

\bibitem{NakatsukasaTropp21}
Y.~Nakatsukasa and J.~A. Tropp.
\newblock Fast \& accurate randomized algorithms for linear systems and
  eigenvalue problems.
\newblock Technical Report arXiv:2111.00113, 2022.

\bibitem{oktay2022mixed}
E.~Oktay and E.~Carson.
\newblock Mixed precision {GMRES}-based iterative refinement with recycling.
\newblock Technical Report arXiv:2201.09827, 2022.

\bibitem{palitta2023sketched}
D.~Palitta, M.~Schweitzer, and V.~Simoncini.
\newblock Sketched and truncated polynomial {K}rylov methods: {E}valuation of
  matrix functions.
\newblock Technical Report arXiv:2306.06481, 2023.

\bibitem{parks2006recycling}
M.~L. Parks, E.~De~Sturler, G.~Mackey, D.~D. Johnson, and S.~Maiti.
\newblock {R}ecycling {K}rylov subspaces for sequences of linear systems.
\newblock {\em SIAM {J}. {S}ci. {C}omput.}, 28(5):1651--1674, 2006.

\bibitem{PhilippeReichel2012}
B.~Philippe and L.~Reichel.
\newblock On the generation of {K}rylov subspace bases.
\newblock {\em Appl.\ Numer.\ Math.}, 62(9):1171--1186, 2012.

\bibitem{rokhlin2008fast}
V.~Rokhlin and M.~Tygert.
\newblock A fast randomized algorithm for overdetermined linear least-squares
  regression.
\newblock {\em Proc. Natl. Acad. Sci. USA}, 105(36):13212--13217, 2008.

\bibitem{Saad1981}
Y.~Saad.
\newblock Krylov subspace methods for solving large unsymmetric linear systems.
\newblock {\em Math.\ Comput.}, 37(155):105--126, 1981.

\bibitem{Saad2003}
Y.~Saad.
\newblock {\em Iterative Methods for Sparse Linear Systems, 2nd edition}.
\newblock SIAM, Philadelphia, 2000.

\bibitem{SaadSchultz1986}
Y.~Saad and M.~Schultz.
\newblock {GMRES:} {A} generalized minimal residual algorithm for solving
  nonsymmetric linear systems.
\newblock {\em SIAM J.\ Sci.\ Stat.\ Comput.}, 7(3):856--869, 1986.

\bibitem{sarlos2006improved}
T.~Sarlos.
\newblock Improved approximation algorithms for large matrices via random
  projections.
\newblock In {\em 47th Annual IEEE Symposium on Foundations of Computer Science
  (FOCS'06)}, pages 143--152. IEEE, 2006.

\bibitem{woodruff2014sketching}
D.~P. Woodruff.
\newblock Sketching as a tool for numerical linear algebra.
\newblock {\em Found. Trends Theor. Comput. Sci.}, 10(1--2):1--157, 2014.

\bibitem{WoolfeLibertyRokhlinTygert2008}
F.~Woolfe, E.~Liberty, V.~Rokhlin, and M.~Tygert.
\newblock A fast randomized algorithm for the approximation of matrices.
\newblock {\em Appl.\ Comput.\ Harmon.\ Anal.}, 25:335--366, 2008.

\end{thebibliography}

\end{document}